\documentclass[a4paper,12pt]{amsart}
\usepackage{amsmath,amsfonts,amssymb}
\usepackage{amsthm}
\usepackage{graphics,graphicx,epsf}
\usepackage[usenames]{color}
\usepackage{color}
\usepackage{hyperref}
\usepackage{tikz}

\newtheorem{theorem}{Theorem}[section]
\newtheorem{definition}[theorem]{Definition}
\newtheorem{proposition}[theorem]{Proposition}
\newtheorem{lemma}[theorem]{Lemma}
\newtheorem{corollary}[theorem]{Corollary}
\newtheorem{remark}[theorem]{Remark}

\renewcommand{\proof}{{\noindent \bf Proof. \ }}

\newcommand{\eproof}{\hfill \mbox{${\square}$}}
\newcommand{\R}{\mathbb{R}}

\newcommand{\sgn}{\textrm{sgn}}

\numberwithin{equation}{section}

\title[Non-autonomous model for neural fields]{Asymptotic behavior for a non-autonomous model of neural fields with variable external stimulus 
}

\author[S. H. da Silva]{Severino Hor\'acio da Silva}
\thanks{Universidade Federal de Campina Grande, Unidade Acad\^emica de Matem\'atica,  58429-900, Campina Grande, PB, Brazil.} 
\thanks{Author E-mail: horacio@mat.ufcg.edu.br and horaciousp@gmail.com - Author Fone: 55-83-21011649.}
\thanks{Partially supported by CNPq/CAPES - Brazil grants Casadinho/Procad \#552.464/2011-2}

\date{\today}

\begin{document}

 \maketitle
\begin{abstract}
In this work we consider a class of nonlocal non-autonomous evolution equations, which generalizes the model of neuronal activity that arises in Amari (1979).


Under suitable assumptions  on the nonlinearity and on the parameters present in the equation, we study, in an appropriated Banach space, the assimptotic behavior of the evolution process generated by this  equation. 
We prove results on existence,  uniqueness and smoothness of the solutions and on the existence of pullback attracts for the evolution process associated to this equation. We also prove a continuous dependence of the evolution process with respect to external stimulus function present in the model.
Furthermore, using the result of continuous dependence of the evolution process, we also prove the upper semicontinuity of pullback attracts with respect to stimulus function. We conclude with a small discussion about the model and about a biological interpretation of the result of continuous dependence of neuronal activity with respect to the external stimulus function.

\vskip .1 in \noindent {\it Mathematical Subject Classification 2010:} 35B40, 35B41, 37B55.\\

\noindent{\it Keywords}: Nonlocal evolution equation; Neural fields; Pullback attractors; Continuous dependence.
\end{abstract}

\section{Introduction}

Neural field equations describe the spatio-temporal evolution of variables such as synaptic or firing rate activity in populations of neurons. The neural field model has already been well  analyzed in the literature 
(see, for example, Amari (1977), 
 Beurle (1956), 
 Bezerra et al. (2012), 
 Carroll and Bressloff (2018), 
 Coombes et al. (2003), 
 Da Silva and Pereira (2009), 
 Da Silva (2010, 2011, 2012), 
 Kishimoto and Amari (1979), 
 Laing (2002), 
 Pinto and Ermentrout (2001), 
 Rukatmatakul and P. Yimprayoon (2008),  Wilson and Cowan (1972),   
 and Zhang (2004). 
  Although this model has been used to model working memory, it arises also in cognitive development of infants, (see, for example, Sandamirskaya (2008) 
  and Thelen (2001)) and 
in  timing sensory integration for robot simulation of autistic behavior (see, for example Barakova (2012)). 

As in Amari (1979), 
we will denote by $u(t,x)$ the membrane potential of neuron located at position $x$ at time $t$, which we are assuming as a differentiable function of $t$, and $J(x,y)$ will denote the average intensity of connections from neurons at place $y$ to those at place $x$.

We also assume that the pulse emission rate of neurons at $x$, at time $t$, is given by a function of $t$ and $u(x,t)$, that is, it is given by  $f(t,u(t,x))$.
 The activity $f(t, u(t,y))$ of neurons at $y$ causes an increase in the potential $u(t,x)$ at $x$ through the connections $J(x,y)$, such that the rate of emission of pulses is proportional to $J(x,y) f(t, u(t,x))$. We also assume that the potential $u(t,x)$ decays, with speed $ 0<\alpha(t) <\alpha_{0}$, to a constant $-h$ (which we call the threshold of the field), while it increases proportionally to the sum of all the stimulus arriving at neurons with speed $b(t)$. Then, denoting by $S(x,t)$ the intensity of the sum of applied stimulus at $x$ at time $t$, and writing $a(t)=\frac{1}{\alpha(t)}$ we have the following non-autonomous evolution equation:

\begin{equation}
\partial_t u(t,x)    =- a(t)u(t,x)  +  b(t) \displaystyle\int_{\mathbb{R}^N} J(x,y)f(t,u(t,y))dy  - h +S(t,x).\label{equation}
\end{equation}
In (\ref{equation}), we consider that   the rate in the intensity of neuronal potential varies explicitly according to the time. Thus, we expect have a more realistic model in (\ref{equation})  when compared to what happens in the brain, since the action potential of electric impulses of the neuronal membrane is consequence of the inversion of polarity inside the membrane, which is not necessarily constant.

Note that, when $a(t)= \frac{1}{\lambda}$, for some constant $\lambda>0$, 
 $b(t)=1$, for all $t\in \mathbb{R}$ and $f(t,x)=f(x)$, equation (\ref{equation}) becomes 
\[
\partial_t u(t,x)    =- u(t,x)  + \displaystyle\int_{\mathbb{R}^N} J(x,y)f(u(t,y))dy -h + S(t,x).
\]
And when $a(t)=b(t)=1$, for all $t\in \mathbb{R}$ and $S(t,x)=h$, equation (\ref{equation}) becomes 
\[
\partial_t u(t,x)    =- u(t,x)  + \displaystyle\int_{\mathbb{R}^N} J(x,y)f(t, u(t,y))dy.
\]
Therefore, equation (\ref{equation}) generalizes the model studied in 
Amari (1977, 1989), 
Bezerra et al. (2017), 
Coombes et al. (2003), 
Da Silva (2010, 2011, 2012), 
Da Silva and Pereira (2015), 
Ermentrout (2010), 
Laing (2002), 
Pinto and Ermentrout (2001), 
Rukatmatakul and P. Yimprayoon (2008), 
Thelen (2001), 
and Zhang (2004). \\

Below we introduce the notations, terminology and some additional  hypotheses, which are already well known in the literature, (see, for example Amari (1977, 1989), 
Bezerra et al. (2017), 
Coombes et al. (2003), 
Da Silva and Pereira (2015), 
and Kishimoto and Amari (1979) ).

Let $\Omega \subset \R^N$ be a bounded smooth domain modelling the geometric configuration of the network,
 $u: \R \times \R^N \to \R$ a function modelling
  the mean membrane potential,  $u(t,x)$ being  the potential
  of a patch of tissue located at position $x\in \Omega$ at time $t\in \R$
 and $f: \R \times \R  \to \R$
 a time dependent transfer function. We say that a neuron  at a point $x$ is active at time $t$  if $f(t,u(t,x))>0$.
 In what follows $b:  \R  \to \R$ is a function in $L^{\infty}([0,\infty))$ such that
$$
0<b(t) \leq b_{0}< \infty,
$$
and it denotes the increasing speed of the potential function $u(t,x)$. Since the decreasing speed of the potential function $u(t,x)$ satisfies $ 0<\alpha(t) <\alpha_{0}$, we can assume that there exist positive constants $a_{-}$ and $a_{0}$ such that  
$$
0<a_{-}\leq a(t) \leq a_0 < \infty.
$$
Let also the integrable function
 $J: \R^N \times \R^N\to \mathbb{R}$  be
 the connection between locations, that is,  $J(x,y)$ is the strength
of the connections of neuronal activity at location $y$ on the activity of the neuron at location $x$. The strength of the connection is assumed to be symmetric,  that is
 $J(x,y)=J(y,x)$, for any $x,y \in \mathbb{R}^N$ and that
 $$
\int_{\mathbb{R}^N} J(x,y) dy = \int_{\mathbb{R}^N} J(x,y) dx= 1.
$$

 Under the conditions assumed above, we analyze  the following   non-autonomous model for neural fields 
\begin{equation} 
\begin{cases}
\partial_t u(t,x)    =- a(t)u(t,x)  + b(t) \displaystyle Kf(t,u(t,y))dy  -h+S(t,x), \ t\geqslant\tau,\  x \in \Omega,\\
u(\tau,x)=u_\tau(x),\ x \in \Omega,\\
 u(t,x) = 0, \ t\geqslant\tau,\  x \in \mathbb{R}^N\backslash\Omega,\label{prob_1}
\end{cases}
\end{equation}
where the integral operator $K$, with symmetric kernel $K$, is given, for all $v\in L ^1(\mathbb{R}^N)$, by
\[
Kv(x):=\int_{{\mathbb{R}^{N}}} J(x,y)v(y)dy.
\]

Also we will assume that $f: \mathbb{R} \times \mathbb{R}\to \mathbb{R}$ satisfies some growth conditions, as presented along  the Section \ref{wellposed}, and that $S: \mathbb{R} \times \mathbb{R}^{N} \to \mathbb{R}$ is continuous in $t$ and $S(t, \cdot) \in L^{p}(\Omega)$, for all $t\in \mathbb{R}$.

\medskip
We aim to study the assimptotic behavior of the evolution process associated to Cauchy problem \eqref{prob_1} in an appropriated Banach space, as well as we obtain some biological conclusion. Then, using the same techniques of Bezerra et al. (2017) 
and Da Silva and Pereira (2015), 
we prove results on existence,  uniqueness and smoothness of the solutions, and the existence of pullback attracts for the evolution process associated to (\ref{prob_1}), which is a more general model than the models analyzed in these previous works. 
We also prove a continuous dependence of the solutions with respect to stimulus function $S$, concluding mathematically that the neuronal activity depends continuously on the sum of external stimulus involved in the neuronal system. This suggests the need for intensive therapies to stimulate people with poor neuronal activity as in some cases of autism or other neurological desorders. Furthermore, using the result of continuous dependence of the evolution process, we also prove the upper semicontinuity of pullback attracts with respect to the stimulus function $S$.

\medskip
This paper is organized as follows. 
In Section \ref{wellposed}, under the growth conditions \eqref{Cf1}, \eqref{Cf2}, \eqref{dissip1} and \eqref{Condf}, for function $f$, we prove  that  \eqref{prob_1} generates a $\mathcal{C}^{1}$ evolution process in the  phase space
\begin{equation}\label{Space}
X_p= \left\{ u \in  L^{p}(\mathbb{R}^{N});\  u(x)= 0,\ \mbox{if}\ x\in \mathbb{R}^N\backslash\Omega \right\}
\end{equation}
with the induced norm, satisfying the ``variation of constants formula''
\[
 u(t,x) =
 \begin{cases}
 e^{-(A(t)-A(\tau))}u_{\tau} (x)+b(t)\displaystyle\int_{\tau}^{t}e^{-(t-s)} K f(s,u(s,\cdot))(x)  ds -h + S(t,x), & x\in \Omega, \\
  0, & x\in \mathbb{R}^N \backslash\Omega,
  \end{cases}
\]
where $A(\xi)=\int_{0}^{\xi} a (\eta) d \eta$, for any $\xi \geq \tau$.
In Section \ref{PullAttractors}, we prove existence of a pullback attractor in the phase space  $X_p$. 
Section \ref{ContParameter} is dedicated to continuity with respect to parameter. In Subsection \ref{ContProcess} we study the continuity of the process with respect to external stimulus function $S$, and we use this result to prove a upper semicontinuity of the pullback attractors in Subsection   \ref{upper-sc}.
Finally, in Section \ref{conclude}, we conclude presenting a brief discussion on the model, with biological interpretation.

\section{The Flow Generated by the Model Problem}\label{wellposed} 

 In this section   we  show the existence of global solution for the problem \eqref{prob_1} and that it generates a $\mathcal{C}^{1}$ evolution processes  in an appropriate Banach space.  For more details on processes evolution  (or infinite-dimensional non-autonomous dynamical systems) see, for example, Carvalho et al. (2012), 
 Chepyzhov and Vishik (2002), 
 Kloeden (2000), 
 and Kloeden and Schmalfu\ss \, (1998). 
  See also Bezerra et al. (2017) 
 and Sell (1967) for related work. 

\subsection{Well posedness}
 In this subsection, under suitable growth condition on the nonlinearity $f$,   we  show the well posedness of the problem \eqref{prob_1} in the  phase space $X_{p}$, for $1\leq p \leq \infty$, given by
 \begin{equation*}
 X_p= \left\{ u \in  L^{p}(\mathbb{R}^{N});\  u(x)= 0,\ \mbox{if}\ x\in \mathbb{R}^N\backslash\Omega \right\}
 \end{equation*}
with the induced norm. It is easy to see that the Banach space $X_p$ is canonically isometric to
 $L^p(\Omega)$, then we usually identify the two spaces, without further comment.  We also use the same notation for a
 function in   $ \R^N$ and its restriction to $\Omega$    for simplicity, wherever we believe the intention is clear from the context.

To obtain well posedness of \eqref{prob_1} in $X_p$, we consider the following Cauchy problem:
\begin{equation}\label{CP}
\begin{cases}
\displaystyle { \frac{d u}{dt}  } =  F(t, u),\ t>\tau, \\
 u(\tau)=u_{\tau},
\end{cases}
\end{equation}
where the map $F: \mathbb{R}\times X_p  \to  X_p $ is defined by
\begin{equation} \label{mapF}
F(t,u)(x)=
\begin{cases}
-a(t)u(x) + b(t) K f(t,u)(x) -h + S(t,x), & \mbox{if}\ t\in\mathbb{R},\ x\in \Omega,\\
 0, & \mbox{if}\ t\in\mathbb{R},\ x\in\mathbb{R}^N\backslash\Omega,
 \end{cases}
\end{equation}
where 
\begin{equation}\label{mapK}
 K f (t, u)(x)  :=   \int_{\mathbb{R}^N} J(x,y)  f (t, u(y)) d y.
\end{equation}
The map $K$ given in (\ref{mapK}) is well defined as a bounded linear operator in  various function spaces, depending on the properties assumed for $J$; for example, with $J$ satisfying the hypotheses from introduction, $K$ is well defined in $X_p$ as shown in lemma below.

\medskip

 The following lemma   has been proved in Da Silva and Pereira (2015).

\begin{lemma} \label{boundK}
Let $K$ be the map defined  by \eqref{mapK} and    $\|J\|_{r}:=\sup_{ x \in \Omega } \|J(x,\cdot)\|_{L^{r} (\Omega) }$, $ 1 \leq r \leq \infty.$ If $u \in   L^p{(\Omega)}, \ 1 \leq p \leq \infty$, then $Ku \in L^{\infty}{(\Omega)}$,
and
\begin{equation}  \label{estimateLq}
|K u (x)|  \leq   \|J\|_{q}      \| u\|_{L^p(\Omega)}\ \mbox{for all}\  x \in \Omega,
\end{equation}
where $1\leq q \leq \infty$ is the conjugate exponent of $p$. Moreover,
\begin{equation}  \label{estimateL1}
   \|K u\|_{L^{p}(\Omega)}   \leq   \|J\|_{1}
   \|  u\|_{L^{p}(\Omega)} \leq
   \|  u\|_{L^{p}(\Omega)}.
   \end{equation}
If $u \in   L^1{(\Omega)} $, then
 $Ku \in L^{p}{(\Omega)}, \ 1 \leq p \leq \infty$,
 and
\begin{equation}  \label{estimateLp}
   \|K u\|_{L^{p}(\Omega)}
 \leq
  \|J\|_{p}
   \|  u\|_{L^{1}(\Omega)}.
   \end{equation}
\end{lemma}

The following definition is already well known in the theory of ODEs in Banach spaces and it can be found in Bezerra at al. (2017). 

\begin{definition} \label{loclip}
If $E$ is a normed space, and $I\subset \R$ is an interval,
we say that a function  $F :I\times E \to E$ is
   \emph{locally Lipschitz
 continuous (or simply locally Lipschitz) } in the second variable  if,
 for any $ (t_0,x_0) \in   I\times E$, there exists  a constant
 $C$  and a rectangle $R = \{ (t,x) \in I \times E \ : \
  |t-t_0|<b_1, \|x-x_0  \|< b_2 \}$ such that,
 if $ (t,x)$ and $(t,y)$ belong to $R$, then
\[
\|F(t,x) - F(t,y) \| \leq  C \| x-y \|.
\]
We say that  $F$ is
   \emph{ Lipschitz
  continuous on bounded sets in the second variable} if the rectangle
  $R$ in the previous definition can chosen as any bounded rectangle in
 $\R \times E$.
 \end{definition}

\begin{remark}
  If the normed space $ E  $ is locally compact the definitions of locally Lipschitz continuous and Lipschitz continuous on bounded sets are equivalent.
\end{remark}

 Now, proceeding as in Bezerra at al. (2017) 
 and Da Silva and Pereira (2015), 
 we prove that the map $F$, given in (\ref{mapF}), is well defined under appropriate growth conditions on $f$ and it is locally  Lipschitz continuous (see Proposition \ref{WellP} below, which generalizes Proposition 3.3 in Bezerra at al. (2017) 
 and Proposition 2.4 in Da Silva and Pereira (2015)).

\medskip
\begin{lemma}\label{WD}
Assume the same hypotheses from Lemma \ref{boundK} and that the function $f$ satisfies the growth condition
	\begin{equation}\label{Cf1}
	|f(t,x)| \leq C_1(t)(1+|x|^{p}),\, \mbox{for any}\ (t,x) \in \mathbb{R}\times\mathbb{R}^N, 
	\end{equation}
	with $1\leq p < \infty$ and $C_1: \mathbb{R} \to \mathbb{R}$ is a locally bounded function.
Then   the function $F$ given
by \eqref{mapF}  is well defined   in $\mathbb{R}\times X_p$. If, for any $t\in \mathbb{R}$, the function $f(t,\cdot)$ is
locally bounded, then $F$ is well defined  in $\mathbb{R}\times L^{\infty}(\Omega)$.
\end{lemma}
\proof
Suppose $1 \leq p < \infty $. Given   $u \in L^p(\Omega)$, denoting the function $f(t,u)(x) = f(t,u(x))$ by  $f(t,u)$ and using \eqref{Cf1}, it easy to see that,  for each $t\in \R$
\begin{equation} \label{estfL1}
\begin{split}
\|f(t,u)  \|_{ L^1(\Omega) } 
& \leq C_1(t)( |\Omega| + \|u\|_{ L^p(\Omega)}^p).
\end{split}
\end{equation}


Thus, using \eqref{estimateLp} and \eqref{estfL1}, it follows that
\begin{eqnarray*}
	\|F(t,u)  \|_{ L^p(\Omega)} &\leq& a_0 \|u\|_{L^p(\Omega)} + b_{0}\|K f(t,u)  \|_{ L^p(\Omega)} + \|S(t, \cdot)\|_{L^{p}(\Omega)} + \|h\|_{L^{p}(\Omega)}\\
	&\leq& a_0 \|u\|_{L^p(\Omega)} +  b_{0} \| J \|_{p}  \|f(t,u)\|_{ L^1(\Omega)} + \|S(t, \cdot)\|_{L^{p}(\Omega)} + h|\Omega|^{\frac{1}{p}}\\
	&\leq& a_0 \|u\|_{L^p(\Omega)} +  b_{0} \| J \|_{p} (C_1(t)|\Omega| + C_1(t) \|u\|_{L^{p}(\Omega)}^{p})  + \|S(t, \cdot)\|_{L^{p}(\Omega)} + h|\Omega|^{\frac{1}{p}}\\
	&\leq& a_0 \|u\|_{L^p(\Omega)} +  b_{0}C_1(t)  \| J \|_{p}  |\Omega| +  b_{0}  C_1(t)  \| J \|_{p} \|u\|_{ L^p(\Omega)}^p  \\
	& +& \|S(t, \cdot)\|_{L^{p}(\Omega)} + h|\Omega|^{\frac{1}{p}}.
\end{eqnarray*}
Since $S(t,\cdot) \in L^{p}(\Omega)$, it follows immediatly that $F$ is  well defined em $L^{p}(\Omega)$ for $1\leq p<\infty$. If $p=\infty$  the result follows easily from  (\ref{estimateLq}).
\qed

\begin{proposition}\label{WellP}
 Under same hypotheses in Lemma \ref{WD}, if $a$ and $b$ are continuous functions and $f$ and $S$, are continuous functions in the first variable, then $F$ is also continuous in the first variable. And if additionally 
\begin{equation}\label{Cf2}
|f(t,x) - f(t,y)| \leq  C_2(t)(1+|x|^{p-1}+|y|^{p-1})|x-y|,\ \mbox{for any}\ (x,y) \in \mathbb{R}^N\times\mathbb{R}^N, \  t\in \mathbb{R},
\end{equation}
for some strictly positive function $C_2:\mathbb{R}\to\mathbb{R}$. Then, for any $1 \leq p < \infty$ the function $F$ is locally Lipschitz continuous on bounded sets in the second variable.
 If $p= \infty$,  this is true if  $f$ is locally Lipschitz in the second variable.
\end{proposition}

\proof
Suppose that $f(t,x)$ is continuous in $t$. Then for any $(t,u) \in \R \times X_p$, we 
 get 
\begin{equation} \label{estfhL1}
\begin{split}
\|f(t,u) -f(t+ \xi, u) \|_{ L^1(\Omega) } & \leq \int_{\Omega}
  |f(t,u(x)) - f(t+ \xi,u(x))|  dx
\end{split}
\end{equation}
for a small $\xi \in \R$.
 From \eqref{Cf1}, follows that the integrand in (\ref{estfhL1}) is bounded by $ 2C (1+|u(x)|^{p}) $, where
 $C$ is a bound for $C(t)$ in a neighborhood of $t$, and goes to
 $0$ as $\xi \to 0$. Hence, using Lebesgue dominated convergence Theorem, it follows that
 $\|f(t,u) -f(t+ \xi, u) \|_{ L^1(\Omega) } \to 0 $ as $\xi \to 0$. Thus, using (\ref{estimateL1}) and (\ref{estfL1}), 
 we obtain
 \begin{eqnarray*}
\|F(t+\xi,u) -F(t,u) \|_{ L^p(\Omega)}&\leq& |a(t)-a(t + \xi)| \|u\|_{L^p(\Omega)}\\ 
&+& |b(t+\xi)-b(t)| \|K (f(t+\xi,u)\|_{ L^p(\Omega)}\\
&+& |b(t)| \|K (f(t+\xi,u)-f(t,u)) \|_{ L^p(\Omega)}\\
&+& \|S(t+ \xi, \cdot) - S(t, \cdot)\|_{L^{p}(\Omega)}\\
&\leq& |a(t)-a(t+\xi)| \|u\|_{L^p(\Omega)} \\
&+&  |b(t+\xi)-b(t)| \| J \|_{p} C_1(t)( |\Omega|+\|u\|_{ L^p(\Omega)}^p)\\
&+& |b(t)| \|J\|_{p} \|f(t+\xi,u)-f(t,u) \|_{ L^1(\Omega)}\\
 &+& \|S(t+ \xi, \cdot) - S(t, \cdot)\|_{L^{p}(\Omega)}
\end{eqnarray*}
which goes to $0$ as $\xi \to 0$, proving the continuity of $F$ in $t$.

Suppose now that
\[
|f(t,x) - f(t,y)| \leq  C_2(t)(1+|x|^{p-1}+|y|^{p-1})|x-y|,
\]
for some $1 < p < \infty$, where $C_2:\mathbb{R}\to\mathbb{R}$ is a strictly positive function. Then, for  $u$ and $v$ belonging to $L^p(\Omega)$ we get
\begin{eqnarray*}
\|f(t,u)-f(t,v)  \|_{ L^{1}(\Omega) } & \leq &  \int_{\Omega} C_2(t)(1+|u(x)|^{p-1} + |v(x)|^{p-1}  )
|u -v | \, d\, x \\
 &  \leq &
  C_2(t) \left[ \int_{\Omega} (1+|u(x)|^{p-1} + |v(x)|^{p-1})^{q} d x \right]^{\frac{1}{q}}
  \left[ \int_{\Omega}   |u(x) - v(x)|^{p} d x  \right]^{\frac{1}{p}}  \\
  &  \leq &
  C_2(t)   \left[ \| 1 \|_{L^q(\Omega)} +
 \| u^{p-1} \|_{L^q(\Omega)}  +  \| v^{p-1} \|_{L^q(\Omega)} \right]
 \|u-v \|_{L^p(\Omega)}
   \\
 &  \leq &
  C_2(t) \left[   | \Omega |^{\frac{1}{q}}  +
 \| u \|_{L^p(\Omega)}^{\frac{p}{q}}  +  \| v \|_{L^p(\Omega)}^{\frac{p}{q}} \right]
   \|u-v \|_{L^p(\Omega)},
\end{eqnarray*}
where $q$ is the conjugate exponent of $p$.

Using  \eqref{estimateLp} once again and the hypothesis on $f$, it follows that
\begin{eqnarray*}
\|F(t,u)- F(t,v)  \|_{ L^p(\Omega)} & \leq & a_0 \|u-v\|_{L^p(\Omega)} + 
  b_0 \|K(f(t,u) - f(t,v))  \|_{ L^p(\Omega)}  \\
 &\leq & a_0 \|u-v\|_{L^p(\Omega)} + b_0 \|J\|_{p} \|f(t,u)-f(t,v)\|_{L^{1}(\Omega)}\\
  &\leq &  \left(a_0 +  b_0 C_2(t) \| J \|_{p} \left[   | \Omega |^{\frac{1}{q}}  +
 \| u \|_{L^p(\Omega)}^{\frac{p}{q}}  +  \| v \|_{L^p(\Omega)}^{\frac{p}{q}} \right]\right)
   \|u-v \|_{L^p(\Omega)},
\end{eqnarray*}
 showing that  $F$  is  Lipschitz in bounded sets of
 ${L^p(\Omega)}$  as claimed.

If $p=1$, the proof is similar. 
Suppose finally that  $\|u\|_{L^{\infty}(\Omega)} \leq R$,
  $\|v\|_{L^{\infty}(\Omega)} \leq R$
 and
 let $M$ be the Lipschitz constant of $f$ in the interval $[-R,R] \subset
 \R $.  Then
 \[
 |f(t,u(x))-f(t,v(x))| \leq M|u(x) - v(x)|,\ \mbox{for any}\ x \in \Omega,
 \]
and this allows us to conclude that
\[
\|f(t,u)-f(t,v)\|_{L^{\infty}(\Omega)}\leq M \|u - v\|_{L^{\infty}(\Omega)}.
\]
Thus,  by \eqref{estimateL1} we have that
  \begin{eqnarray*}
  \|F(t,u)-F(t,v)\|_{L^{\infty}(\Omega)} & \leq & a_0 \|u-v\|_{L^{\infty}(\Omega)} +
b_0 \|K (f(t,u) - f(t,v)) \|_{L^{\infty}(\Omega)}\\
   &  \leq & \left(a_0 + b_0 M \|J\|_1 \right)  \|u -v \|_{L^{\infty}(\Omega)},
  \end{eqnarray*}
and this completes the proof.
\qed

\bigskip


Using Proposition \ref{WellP} and well known results of  ODEs in Banach spaces, see Daleckii and Krein (1974), 
 it follows  that  the initial value problem \eqref{CP} has a unique
local solution for any initial condition in $X_p$.
For the  global existence, we need of the Theorem 5.6.1 in Ladas and Lakshmikantham (1972). 


\begin{proposition} \label{globalexist}
	Besides the assumptions from Proposition \ref{WellP} we suppose also that there exists, constant $k_1 \in \R$, independent of $t$, such that 	
	$ f $  satisfies the dissipative condition
\begin{equation}\label{dissip1}
 \limsup_{|x| \to \infty} \frac{|f(t,x)|}{|x|} < k_1.
 \end{equation}
 Then the   problem  (\ref{CP})  has a unique   globally defined
 solution for any initial
 condition in  $X_p$,
 which is  given, for $t\geq \tau$, by
 the ``variation of constants formula''

\begin{equation}\label{EP_1}
u(t,x) =
\begin{cases}
e^{-(A(t) - A(\tau))}u_{\tau} (x)+\displaystyle\int_{\tau}^{t}e^{-(A(t) - A(s))}[b(s) K f(s,u(s,\cdot))(x) -h + S(s,x)] ds,\!\!\!\! 
&x\in \Omega, \\ 
0,\!\!\!\!\!\!\!\!\!\!\!\!&x\in \Omega^c,
\end{cases}
\end{equation}
where $A(\xi)=\int_0^\xi a(\eta) d \eta$, for any $\xi \geq \tau$, and $ \Omega^c =\mathbb{R}^N \backslash\Omega$.

\end{proposition}

\proof
Existence and uniqueness of local solutions for (\ref{CP}), in $X_p$, it follows from Proposition \ref{WellP} and well-known results (see Daleckii and Krein (1974)). 
The variation of constants formula (\ref{EP_1}) can be easily  verified by direct derivation.


Now, using condition \eqref{dissip1} it follows that
\begin{equation} \label{dissip2}
| f(t,x)  | \leq k_2(t) + k_1 |x |,\ \mbox{for any}\ (t,x) \in \mathbb{R}\times\mathbb{R}^N,
\end{equation}
for some continuous and strictly positive function $k_2:\mathbb{R}\to\mathbb{R}$.

If $ 1 \leq p < \infty$,  using \eqref{estimateL1}   and \eqref{dissip2}, we obtain  the following estimate
 \begin{eqnarray*}
  \| K f(t, u )\|_{L^p(\Omega)} &   \leq &   \|f(t, u)\|_{L^p(\Omega)}\\
   &   \leq & k_2(t) |\Omega|^{1/p} + k_1 \|   u \|_{L^p(\Omega)}.
\end{eqnarray*}

For $p= \infty$, using the  same arguments (or by making $p \to \infty$ in previous inequality), we have
\[
  \| K f(t,u )\|_{L^{\infty}(\Omega)} \leq   k_2(t)  + k_1 \|    u  \|_{L^{\infty}(\Omega)}.
\]

Now defining the function
\[
\begin{split}
g: [t_0, \infty) \times  \R^{+}  &\to   \R^{+}\\
(t,r) &\mapsto g(t,r)= b_{0}  |\Omega|^{1/p}  k_2(t) + \|S\|_{p} + h |\Omega|^{1/p}   +
 ( k_1 +a_0) r
\end{split}
\]
it follows that problem \eqref{CP} satisfies the hypothesis of  
Theorem 5.6.1 in Ladas and Lakshmikantham (1972) 
and the global existence
 follows immediately. 

\subsection{Smoothness of the evolution process} \label{Smoothness}

 In this subsection   we  show that the problem \eqref{prob_1} generates a $\mathcal{C}^{1}$ flow in the  phase space $X_{p}$.\\

%

\begin{proposition} \label{C1flow}
Assume the same hypothese from Proposition \ref{globalexist} and that 
 the function $f$ is  continuously differentiable in the second variable
and $\partial_2f$ satisfies the growth condition
\begin{equation}\label{Condf}
|\partial_2f(t,x)| \leq C_1 (t)(1+|x|^{p-1}),\ \mbox{for any}\ (t,x) \in \mathbb{R}\times\mathbb{R}^N,
\end{equation}
for $1  \leq p < \infty$. Then $F(t, \cdot)$ is continuously Frech\'et differentiable on
$X_p$ with derivative given by
\begin{eqnarray*}
	DF(t,u)v(x):=
	\begin{cases}
	-a(t)v(x) +	b(t)K(\partial_2f (t,u)v)(x), & x\in \Omega,\\
		0, & x\in \mathbb{R}^N\backslash\Omega.
 \end{cases}
\end{eqnarray*}
\end{proposition}
\proof Using that $f$ is continuously differentiable in the second variable, by a simple computation, it follows  that the Gateaux's derivative of $F(t, \cdot)$ is given by
\begin{eqnarray*}
DF(t,u)v(x):=
\begin{cases}
	-a(t)v(x) + b(t)K(\partial_2f(t,u)v)(x)
, & x\in \Omega,\\
 0, & x\in \mathbb{R}^N\backslash\Omega,
 \end{cases}
\end{eqnarray*}
where  $(\partial_2f(t,u)v)(x) := \partial_2f(t,u(x))\cdot v(x)$.
Note that the operator $D_2F(t,u)$
 is a linear operator in $X_p$.

Let $u \in L^p(\Omega)$, with $1 \leq p < \infty $. Then, if $q$ is the conjugate exponent of $p$, it is easy to see that

\begin{eqnarray} \label{estflinhaq}
\|\partial_2f(t,u)  \|_{ L^q(\Omega) } \leq   C_1(t) \left( |\Omega|^{\frac{1}{q}} +  \|u \|_{ L^p(\Omega)}^{p-1} \right).
\end{eqnarray}
From H\"{o}lder inequality and estimate (\ref{estflinhaq}), it follows that 
\[
\|\partial_2f(t,u)\cdot v   \|_{ L^1(\Omega) } \leq
  C_1(t)( |\Omega|^{\frac{1}{q}} +  \|u \|_{ L^p(\Omega)}^{p-1} ) \|v\|_{ L^p(\Omega)}.
\]

Hence, from estimate  \eqref{estimateLp}, we concluded that
\begin{eqnarray*}
\|DF(t,u)\cdot v  \|_{ L^p(\Omega)} &\leq& a_0 v + b_{0}\|K(\partial_2f(t,u)v)  \|_{ L^p(\Omega)}\\
& \leq& a_0 v +  b_{0} C_1(t) \| J \|_{p}  \|\partial_2f(t,u)v\|_{ L^1(\Omega)}\\
 &\leq& a_0 v +  b_{0} C_1(t) \| J \|_{p}
 \left( |\Omega|^{\frac{1}{q}} +  \|u \|_{ L^p(\Omega)}^{p-1}  \right)\|v\|_{ L^p(\Omega)}\\
&=& [a_0  +  b_{0} C_1(t) \| J \|_{p}
\left( |\Omega|^{\frac{1}{q}} +  \|u \|_{ L^p(\Omega)}^{p-1} \right)]\|v\|_{ L^p(\Omega)} ,
\end{eqnarray*}
that is, $DF(t,u)$ is
a bounded operator.
In the case $p= \infty$, it follows that $ |\partial_2f(t,u)|$ is bounded by  $C_2(t)$,
 for each $ u  \in L^{\infty}(\Omega)$.
 Hence
 \[
 \|\partial_2f(t,u) v \|_{ L^{\infty}(\Omega)} \leq C_2(t) \| v \|_{ L^{\infty}(\Omega)}.
 \]
 Thus, using \eqref{estimateL1}, we obtain
 \[
\begin{split}
\|DF(t,u)\cdot v  \|_{ L^{\infty}(\Omega)} &\leq a_0 \|v\|_{L^{\infty}} + b_{0} \|K(\partial_2f(t,u)v)  \|_{ L^{\infty}(\Omega)}\\
& \leq a_0 \|v\|_{L^{\infty}} +  b_{0}  \| J \|_{1}  \|\partial_2f(t,u)v\|_{ L^{\infty}(\Omega)}\\
& \leq a_0 \|v\|_{L^{\infty}} + b_{0}  C_2 (t) \| J \|_{1}  \| v \|_{ L^{\infty}(\Omega)}\\
& = (a_0  + b_{0}  C_2 (t) \| J \|_{1} ) \| v \|_{ L^{\infty}(\Omega)},
\end{split}
\]
 which results in the boundedness of  $DF(t,u)$ also in this case.

 Now, suppose that $u_1$ and $u_2$ and $v$  belong to $L^p(\Omega)$,
 $ 1 \leq p < \infty$.
Using \eqref{estimateLp}  and H\"{o}lder inequality, it follows  that
 \begin{eqnarray*}
 \| (DF(t,u_1) - DF(t,u_2)) v\|_{L^p(\Omega)} & \leq & 
 b_{0} \|K [ (\partial_2f(t,u_1) - \partial_2f(t,u_2)) v    ]   \|_{L^p(\Omega)} \\
 & \leq & 
 b_{0} \| J  \|_p  \|  (\partial_2f(t,u_1) - \partial_2f(t,u_2)) v       \|_{L^1(\Omega)} \\
 & \leq &  
 b_{0} \| J  \|_p  \|  \partial_2f(t,u_1) - \partial_2f(t,u_2)     \|_{L^q(\Omega)}  \| v \|_{L^p(\Omega)}\\
 & = &  
 b_{0} \| J  \|_p  \|  \partial_2f(t,u_1) - \partial_2f(t,u_2)     \|_{L^q(\Omega)}  \| v \|_{L^p(\Omega)}.
  \end{eqnarray*}
Then to prove continuity of the derivative, $DF(t,\cdot)$, it is sufficient to show that
\[
\|  \partial_2f(t,u_1) - \partial_2f(t,u_2) \|_{L^q(\Omega)} \to 0
\]
as $ \|  u_1 - u_2     \|_{L^p(\Omega)} \to 0 $.
But, from  (\ref{Condf}), it follows that
\[
|\partial_2f(t,u_1)(x)- \partial_2f(t,u_2)(x)|^q \leq [C_1(t)(2+|u_1(x)|^{p-1} + |u_2(x)|^{p-1})]^q.
\]
But a simple computation shows  
that the right-hand-side of this last inequality is integrable.  Then the result follows from Lebesgue's convergence
 theorem.

In the case $p =\infty$, the continuity of $DF$ follows from \eqref{estimateL1} and from the continuity of $\partial_2f(t,u)$.

Therefore, it  follows from  Proposition 2.8 in Rall (1971) 
that $F(t, \cdot)
$ is Fr\'echet
differentiable  with  continuous  derivative in  $X_p$.\qed\\

As a consequence of the Proposition \ref{C1flow} and of well know results of Daleckii and Krein  (1974) 
and Henry  (1981), 
we have the following result.
\begin{corollary} \label{Smoothness}
Assume the same hypotheses of the Proposition \ref{C1flow}. Then,  
for each $t\in \mathbb{R}$ and $u_\tau \in X_p$, the unique solution of \eqref{CP} with initial condition $u_\tau$ exists for all $t\geq \tau$ and   this solution $(t,\tau, x)\mapsto u(t,x)=u(t;\tau,x,u_\tau)$ (defined by \eqref{EP_1}) gives rise to a family of nonlinear $\mathcal{C}^1$ process on $X_p$, given by
\[
T(t,\tau)u_{\tau}(x):=u(t,x),\ t\geq \tau\in \R.
\]
\end{corollary}

\section{Existence  of pullback attractor }   \label{PullAttractors}

In this section we   prove the existence of a pullback attractor $\{\mathcal{A}(t); t\in\mathbb{R}\}$ in $X_p$ for the evolution process $ \{T(t, \tau); t \geq \tau, \tau\in\mathbb{R}\} $ when $ 1 \leq p  <  \infty $, generalizing, among others, Theorem 3.2 of Da Silva and Pereira (2015) 
and Theorem 4.2 of Bezerra at al. (2017). 

\begin{lemma}\label{L_PullAbs}
Assume that the hypotheses from Proposition \ref{C1flow} hold with the constant
 $k_1$ in (\ref{dissip1}) satisfying  $k_{1}b_0 < a_{-}$. Let
 \begin{equation}\label{DefR0}
 R_\delta(t)=\dfrac{1}{a_{-}-k_1 b_0}(1+\delta) [b_0 k_2(t) |\Omega|^{\frac{1}{p}} +\| S(t,\cdot)\|_{L^p(\Omega)}],
 \end{equation}
  where $k_2$ is derived from \eqref{dissip2} and $\delta$ is any positive constant.
  Then the ball of $L^p(\Omega), \ 1 \leq p <  \infty$, centered
at the origin with  radius $R_\delta(t)$,
which we denote by $\mathcal{B}(0,R_\delta(t))$, pullback  absorbs bounded subsets of $X_p$ at time $t\in\mathbb{R}$ with respect to the  process $T(\cdot,\cdot)$ generated by \eqref{CP}.
\end{lemma}
 \proof  If  $u(t,x)$ is the solution of \eqref{CP} with
initial condition $u_{\tau} \in X_p$, for $ 1 \leq p <\infty$, then 
 \begin{equation}
\begin{split}
&\frac{d}{dt}\int_{\Omega}|u(t,x)|^{p}dx
=  \int_{\Omega} p |u(t,x)|^{p-1} \sgn(u(t,x))
   u_t(t,x) dx  \\
&= -p a(t)\int_{\Omega}|u(t,x)|^{p}dx +
p b(t) \int_{\Omega}|u(t,x)|^{p-1} \sgn(u(t,x))K f(t,u(t,x)) dx\\
&+ p \int_{\Omega}|u(t,x)|^{p-1} \sgn(u(t,x)) S(t,x) dx - p h \int_{\Omega}|u(t,x)|^{p-1}  dx.
\end{split}     \label{mkq1}
\end{equation}

Thus, if $q$ is the conjugate exponent of $p$, from H\"{o}lder's inequality, estimate \eqref{estimateL1} and condition \eqref{dissip1},  we have
 \begin{equation}\label{mkq2}
\begin{split}
&\int_{\Omega}|u(t,x)|^{p-1} \sgn(u(t,x)) K f(t,u(t,x))dx\\
 & \leq  \left(\int_{\Omega} |u(t,x)|^{q(p-1)} dx\right)^{\frac{1}{q}}
 \left(
\int_{\Omega} | K f(t,u(t,x))|^pdx \right)^{\frac{1}{p}}  \\
 & \leq  \left(\int_{\Omega} |u(t,x)|^{p}dx\right)^{\frac{1}{q}}
 \|J\|_1 \| f(t,u(t,\cdot))  \|_{L^{p}(\Omega)}
   \\
 &\leq  \|u(t,\cdot)\|_{L^{p}(\Omega)}^{p-1}
 \left(k_1 \|u(t,\cdot)\|_{L^{p}(\Omega)} + k_2(t) |\Omega|^{\frac{1}{p}}\right),
\end{split}
\end{equation}
and
\begin{equation}\label{mkq3}
\begin{split}
&\int_{\Omega}|u(t,x)|^{p-1} \sgn(u(t,x)) S(t,x)dx\\
 & \leq  \left(\int_{\Omega} |u(t,x)|^{q(p-1)} dx\right)^{\frac{1}{q}}
 \left(
\int_{\Omega} | S(t,x)|^pdx \right)^{\frac{1}{p}}  \\
 & \leq  \left(\int_{\Omega} |u(t,x)|^{p}dx\right)^{\frac{1}{q}}
\|S(t,\cdot)  \|_{L^{p}(\Omega)}
   \\
 &\leq  \|u(t,\cdot)\|_{L^{p}(\Omega)}^{p-1}
 \|S(t,\cdot)\|_{L^{p}(\Omega)}.
\end{split}
\end{equation}
Hence, using (\ref{mkq2}) and (\ref{mkq3}) in (\ref{mkq1}), we obtain

\[
\begin{split}
\frac{d}{dt}\|u(t,\cdot)\|_{L^{p}(\Omega)}^{p} &\leq -p a(t) \|u(t,\cdot)\|_{L^{p}(\Omega)}^{p}+p b(t) \|u(t,\cdot)\|_{L^{p}(\Omega)}^{p-1} \big( k_1 \|u(t, \cdot) \|_{L^{p}(\Omega)} + k_2(t){|\Omega|}^{\frac{1}{p}}  \big)\\
& + p\|u(t,\cdot)\|_{L^{p}(\Omega)}^{p-1} \|S(t,\cdot)\|_{L^{p}(\Omega)} - p h |\Omega|^{\frac{1}{p}} \|u(t,\cdot)\|_{L^{p}(\Omega)}^{p-1}.
\end{split}
\]
Thus
 \[
\begin{split}
\frac{d}{dt}\|u(t,\cdot)\|_{L^{p}(\Omega)}^{p} &\leq
-a_{-}p \|u(t,\cdot)\|_{L^{p}(\Omega)}^{p} + p b_{0} \|u(t,\cdot)\|_{L^{p}(\Omega)}^{p-1} \left(k_{1} \|u(t,\cdot)\|_{L^{p}(\Omega)} +
k_{2}(t) |\Omega|^{\frac{1}{p}}\right)\\
& + p  \|S(t,\cdot)\|_{L^{p}(\Omega)} \|u(t,\cdot)\|_{L^{p}(\Omega)}^{p-1} \\
&=-a_{-}p \|u(t,\cdot)\|_{L^{p}(\Omega)}^{p} + p b_{0} k_1 \|u(t,\cdot)\|_{L^{p}(\Omega)}^{p} +  p b_{0} |\Omega|^{\frac{1}{p}}  k_{2}(t) \|u(t,\cdot)\|_{L^{p}(\Omega)}^{p-1}  \\
& + p  \|S(t,\cdot)\|_{L^{p}(\Omega)} \|u(t,\cdot)\|_{L^{p}(\Omega)}^{p-1}\\
&=p \|u(t,\cdot)\|_{L^{p}(\Omega)}^{p} \left[-a_{-} + k_{1}b_{0} + \frac{( b_{0} k_{2}(t) |\Omega|^{\frac{1}{p}} +  \|S(t,\cdot)\|_{L^{p}(\Omega)} )}{ \|u(t,\cdot)\|_{L^{p}(\Omega)}^{p} }\right].
\end{split}
\]
Writing $\varepsilon =a_{-}-k_{1}b_{0} > 0$, it follows that, 
while
\[
\|u(t,\cdot)\|_{L^{p}(\Omega)} \geq \dfrac{1}{\varepsilon}(1+\delta) \left( b_{0} k_{2}(t) |\Omega|^{\frac{1}{p}} +   \|S(t,\cdot)\|_{L^{p}(\Omega)} \right),
\]
we obtain
\begin{eqnarray} \label{boundsol}
\frac{d}{dt}\|u(t,\cdot)\|_{L^{p}(\Omega)}^{p} &\leq
p \|u(t,\cdot)\|_{L^{p}(\Omega)}^{p} \left(-\varepsilon + \frac{\varepsilon}{1+ \delta} \right) \nonumber \\
&= -\frac{\delta p}{(1+\delta)} \varepsilon \| u(t, \cdot)\|_{L^{p}(\Omega)}^{p}.\nonumber
\end{eqnarray}

Therefore
\begin{eqnarray} \label{boundsol}
\|u(t,\cdot)\|_{L^{p}(\Omega)}^{p} &\leq&
e^{-\frac{ \delta  p  }{(1+\delta)}  \varepsilon (t-\tau)} \| u_{\tau}\|_{L^{p}(\Omega)}^{p} \nonumber \\
&=&e^{- \frac{\delta p}{(1+\delta)}(a_{-} - k_1 b_{0})(t-\tau)}\| u_{\tau}\|_{L^{p}(\Omega)}^{p}.
\end{eqnarray}
Thus, the  result follows immediately.
 \qed

\begin{theorem}\label{Theor1}
In addition to the  conditions of  Lemma \ref{L_PullAbs}, suppose that $C_1 (t)$ and $K_2(t)$ are non-decreasing functions and 
\[
\| J_x \|_{L^p(\Omega)} = \sup_{x\in\Omega}\|\partial_x J(x,\cdot)\|_{L^q(\Omega)}< \infty; \,\,
\mbox{and} \,\, 
\|\partial_{x}S \|_{p}=   \sup_{t\in \mathbb{R}_{+}} \|\partial_{x}S(t, \cdot)\|_{L^{p}(\Omega)} < \infty.
\]
Then there exists a pullback attractor $\{\mathcal{A}(t); t\in\mathbb{R}\}$ for the process $\{T(t,\tau); t\geq \tau, \tau\in\mathbb{R}\}$ generated by \eqref{CP} in $X_p= L^p(\Omega)$  and the `section' $\mathcal{A}(t)$ of the pullback attractor $\mathcal{A}(\cdot)$ of $T(\cdot,\cdot)$ is contained in the ball centered at the origin with  radius $R_\delta(t)$ defined in \eqref{DefR0}, in $L^p(\Omega)$,  for any $\delta>0$, $t \in \mathbb{R}$ and $1 \leq p < \infty$.
\end{theorem}

\proof
From Theorem \ref{globalexist} it follows that, for each initial value $u(\tau,\cdot)\in X$ and initial time $\tau\in\mathbb{R}$, the process generated by \eqref{CP} has a unique solution, which we can to write, for $x\in \Omega$, as
\[
T(t,\tau)u(\tau,x) =T_{1}(t,\tau) u(\tau,x)  + T_{2}(t,\tau) u(\tau,x),
\]
where 
\[
T_{1}(t,\tau)u(\tau,x):=e^{-(A(t)-A(\tau))} u(\tau,x),
\]
 and
\[
T_{2}(t,\tau)u(\tau,x):=\int_{\tau}^{t}e^{-(A(t)-A(s))} b(s) [ Kf(s, u(s,x))  + S(s,x) -h ]ds.
\]

Now, we will use the Theorem 2.37 in Carvalho at al. (2012), 
to prove that $T(\cdot,\cdot)$ is pullback asymptotically compact. For this, let $u\in B$ be, where $B$ is a bounded subset of $X_p$. Without loss of generality, we may suppose that $B$ is contained in the ball  centered at the origin of radius $ r>0$. Then, for $t\geq \tau$, we have
\begin{eqnarray}
\|T_{1}(t,\tau) u\|_{L^p(\Omega)} &\leq& 
 r e^{-(A(t)-A(\tau))}  \nonumber \\
&\leq& r e^{-a_{-}t} e^{a_{0}\tau}=\sigma(t,\tau) \to 0, \, t \to \infty  \nonumber .
\end{eqnarray}
Using (\ref{boundsol}), it follows that $\|u(t,\cdot)\|_{L^p(\Omega)}\leq M$, for $t\geq \tau$, where $M$ is given in (\ref{est_attractor}) below
\begin{equation}
M=M(t)= \max \left\{r,  \frac{2 [b_0 k_2(t)|\Omega|^{\frac{1}{p}} + \|S(t, \cdot) \|_{L^{p}(\Omega)}] }{a_{-}-k_1 b_0} \right\}>0. \label{est_attractor}
\end{equation}

Then, using  (\ref{estfL1}), we have
\[
\begin{split}
\|f(t,u)  \|_{ L^1(\Omega) } &\leq C_1(t)( |\Omega| + \|u\|_{ L^p(\Omega)}^p)\\
& \leq C_1(t)( |\Omega| + M(t)^p).
\end{split}
\]
Since 
\[
\partial_{ x} ( T_2 (t,\tau) u(\tau, x))=\int_{\tau}^{t} e^{-(A(t)-A(s))} [ b(s) \frac{\partial }{\partial x} Kf(t,u)(t,x) + \frac{\partial S}{\partial x} (s,x) ] ds.
\]
proceeding as in  \eqref{estimateLp} (with $J_x$ in the place of $J$) and using (\ref{estfL1}), it follows that
\[
\begin{split}
\left\| \partial_{ x}( K f(t,u)) \right\|_{L^p(\Omega)}
 &\leq  \| J_x \|_{L^p(\Omega)}  b_{0} \|f(t, u)\|_{ L^1(\Omega)}\\
& \leq  C_1(t) \| J_x \|_{L^p(\Omega)} ( |\Omega| + M(t)^p ).
\end{split}
\]
Thus, since $C_1$ and $k_2$ are non-decreasing, we obtain
\begin{eqnarray*}
\label{est_deriv_T_2}
\left\| \partial_{ x} (T_2(t,\tau)u ) \right\|_{L^{p}(\Omega) }
& \leq& \int_{\tau}^{t}e^{-(A(t)-A(s))} \left( b(s)\| \partial_{x}K f(s,u(s,\cdot))\|_{L^{p}(\Omega) }  + \left\| \partial_{ x} S(s,\cdot) \right\|_{L^{p}(\Omega) } \right) ds\\
& \leq&
\int_{\tau}^{t}e^{-(A(t)-A(s))}  \left( b_0 C_1(s)\| J_{x}\|_{L^{p}(\Omega) }  |\Omega|^{\frac{1}{p}} + M(s)^{p} + \left\| \partial_{ x} S \right\|_{p } \right)ds\\
& \leq&
\int_{\tau}^{t}e^{-(A(t)-A(s))}  \left( b_0 C_1(t)\| J_{x}\|_{L^{p}(\Omega) }  |\Omega|^{\frac{1}{p}} + M(t)^{p} + \left\| \partial_{ x} S \right\|_{p } \right)ds\\
& \leq&   C_1(t)\| J_x \|_{p}  \frac{1}{a_0} \left[ e^{(a_0 -a_{-})t} -  e^{-a_{-}t}e^{a_{0} \tau} \right]  \left( |\Omega|^{\frac{1}{p}} + M(t)^p \right)\\
&+& C_1(t)  \frac{1}{a_0} \left[ e^{(a_0 -a_{-})t} -  e^{-a_{-}t}e^{a_{0} \tau} \right]  \left\| \partial_{ x} S \right\|_{p }\\
& \leq &  C_1(t)\| J_x \|_{p}  \frac{1}{a_0}  e^{(a_0 -a_{-})t}   \left( |\Omega| + M(t)^p \right) + C_1(t)  \frac{1}{a_0}  e^{(a_0 -a_{-})t}   \left\| \partial_{ x} S \right\|_{p }\\
&=&\frac{ C_1(t)\| J_x \|_{p} ( |\Omega| + M(t)^p ) + \left\| \partial_{ x} S \right\|_{p } }{a_0} e^{(a_0 -a_{-})t} .
\end{eqnarray*}
Hence, for any $ u \in B$ and $t>\tau$, the value of
$\|\frac{\partial }{\partial x} T_2(t,\tau)u \|_{L^{p}(\Omega)}$ is bounded by
a constant (independent of $u \in B$). Then  $T_2(t,\tau) u $ belongs to  a ball of $W^{1,p}(\Omega)$  for all $u \in B$. Hence,
 from Sobolev's Embedding Theorem, it  follows that
  $T_2(t,\tau) $ is a compact operator, for any $ t>\tau$.

Therefore, using Lemma \ref{L_PullAbs} and Theorem 2.23 of Carvalho at al. (2012), 
it follows that there exists the
pullback attractor $\{\mathcal{A}(t); t\in\mathbb{R}\}$ and each `section' $\mathcal{A}(t)$ of the pullback attractor $\mathcal{A}(\cdot)$ is the pullback $\omega$-limit set of any bounded subset of $X_p$ containing the ball centered at the origin with  radius $R_\delta$, given in \eqref{DefR0}, for any $\delta > 0$. Since  the ball centered at the origin with  radius $R_\delta$ pullback absorbs bounded subsets of   $X_p$,  it also follows that the set $ \mathcal{A}(t)$  is contained in the ball centered at the origin of $X_p$ and of radius
\[
R(t)=\dfrac{1}{a_{-}-k_1 b_0} \left[b_0 k_2(t) |\Omega|^{\frac{1}{p}} +\left\|  S \right\|_{p}\right], \,\, \mbox{ for any} \,\, t \in \mathbb{R}, \,\, 1 \leq p < \infty.
\]
\qed

\vspace{0.7cm}
\section{Continuity with respect to parameter} \label{ContParameter}
\vspace{0.5cm}
A natural question to examine is the depedence of the process with respect to parameters that arise in the equation.
In this section we prove the continuity of the process with respect to external stimulus function and we use this result to prove the upper semicontinuity of the pullback attractors.
\vspace{0.6cm}
\subsection{Continuity of the process with respect to external stimulus}\label{ContProcess}

From now on we denote by $T_{S}(t,\tau)$ the family of processes associated with the family of problems

\begin{equation}\label{pn1}
\begin{cases}
\partial_t u_S(t,x)    =- a(t)u_S(t,x)  + b(t) K f(t,u_{S}(t,x)) + S(t,x), \ t\geq\tau,\  x \in \Omega,\\
u_S(\tau,x)=u_{\tau}(x),\ x \in \Omega,\\
u_S(t,x) = 0, \ t\geq\tau,\  x \in \mathbb{R}^N\backslash\Omega.
\end{cases}
\end{equation}

In this subsection we prove the continuous dependence of the process with respect stimulus function $S$ at $S_{0} \in \Sigma$, where
$$
\Sigma=\{ S: \mathbb{R} \times \mathbb{R}^N \to \mathbb{R}, \, \|S\|_{p}=\sup_{t\in \mathbb{R}_{+}}\| S(t,\cdot)\|_{L^{p}(\Omega)} < \infty
\}.
$$
More precisely we have the following result:

\begin{theorem}\label{tend}
	In addition to the hypotheses of  Theorem \ref{Theor1}, suppose that the function $C_2 (t)$ given in (\ref{Cf2}), is non-decreasing. Then, if $T_S(\cdot,\cdot)$ denotes the process generates by the problem \eqref{pn1}, for $S\in \Sigma$, we have that	
	\[
	\mbox{$T_S(t,\tau)u_{\tau}$ converges to  $T_{S_0}(t,\tau)u_{\tau}$ in $X_p$, as $\|S-S_0\|_{p} \to 0$},
	\]
	uniformly for $t\in [\tau,L]$, for any $L>\tau$.
\end{theorem}

\proof
Let $L>\tau $ and $u_S(t,x)=T_S(t,\tau)u_{\tau}(x)$ be the solution of the problem \eqref{pn1} for $t\in [\tau,L]$, given by \eqref{EP_1}. Then, for $x\in \Omega$,
\begin{eqnarray*}
u_S(t,x)-u_{S_0}(t,x)&=&\int_{\tau}^{t}e^{-(A(t)-A(s))} b(s) [K( f(s,u_S(s,x)) - f(s,u_{S_0}(s,x)))] ds\\
&+& \int_{\tau}^{t}e^{-(A(t)-A(s))} [S(s,x) -S_{0}(s,x)]ds
\end{eqnarray*}
Thus, for $x\in \Omega$, using  \eqref{estimateLp}, we obtain 
\[
\begin{split}
\|u_S(t,\cdot)-u_{S_0}(t,\cdot)\|_{L^p(\Omega)} &
\leq\displaystyle\int_{\tau}^{t}e^{-(A(t)-A(s))}  b_0  \|J\|_{p}\|  f(s,u_S(s,\cdot))- f(s,u_{S_0}(s,\cdot))\|_{L^1(\Omega)}ds\\
&+\displaystyle\int_{\tau}^{t}e^{-(A(t)-A(s))}    \| S(s,\cdot) -S_{0}(s,\cdot)\|_{L^p(\Omega)}ds.
\end{split}
\]
From \eqref{estfL1} it follows that
\[
\begin{split}
\|u_S(t,\cdot)-u_{S_0}(t,\cdot)\|_{L^p(\Omega)} &
\leq\displaystyle\int_{\tau}^{t}e^{-(A(t)-A(s))}  b_0  \|J\|_{p} C_{2} (s) \bigg[|\Omega|^{\frac{1}{q}} + \|u_{S}(s,\cdot)\|_{L^{p}(\Omega)}^{\frac{p}{q}}\\ 
&+ \|u_{S_0}(s,\cdot)\|_{L^{p}(\Omega)}^{\frac{p}{q}} \bigg]\| u_S(s,\cdot) - u_{S_0}(s,\cdot)\|_{L^p(\Omega)}ds\\
&+\displaystyle\int_{\tau}^{t}e^{-(A(t)-A(s))}   \sup_{s\in \mathbb{R}} \| S(s,\cdot) -S_{0}(s,\cdot)\|_{L^p(\Omega)}ds.
\end{split}
\]

Let  $B \subset X_p$ a bounded subset  (for example a ball of radius $\rho$) such that $u_S(t,\cdot)\in B$  for all $S\in \Sigma$ and $t\in[\tau,L]$. 
Then, 
we have
\[
\begin{split}
e^{A(t)}\|u_S(t,\cdot)-u_{S_0}(t,\cdot)\|_{L^p(\Omega)} &
\leq\displaystyle\int_{\tau}^{t}  b_0  \|J\|_{p} C_2(s) \bigg[|\Omega|^{\frac{1}{q}} + 2\rho^{\frac{p}{q}} \bigg]e^{A(s)}\| u_S(s,\cdot) - u_{S_0}(s,\cdot)\|_{L^p(\Omega)}ds\\
&+\displaystyle\int_{\tau}^{t}e^{A(s)}  \| S -S_{0}\|_{p}ds.
\end{split}
\]
From Generalized Gronwall inequality, see Hale (1980), we get
\[
\begin{split}
e^{A(t)}\|u_S(t,\cdot)-u_{S_0}(t,\cdot)\|_{L^p(\Omega)} &
\leq\displaystyle\left(\int_{\tau}^{t} e^{A(s)} \| S -S_{0}\|_{p}ds \right) e^{ \int_{\tau}^{t}b_0  \|J\|_{p} C_2(s) \left[ |\Omega|^{\frac{1}{q}} + 2\rho^{\frac{p}{q}}\right] ds}.
\end{split}
\]
Hence, for $t\in[\tau,L]$, it follows that

\[
\begin{split}
\|u_S(t,\cdot)-u_{S_0}(t,\cdot)\|_{L^p(\Omega)} &
\leq\displaystyle\left(\int_{\tau}^{t} e^{-(A(t)-A(s))} \| S -S_{0}\|_{p}ds \right) e^{ \int_{\tau}^{t} b_0  \|J\|_{p} C_2(s) \left[ |\Omega|^{\frac{1}{q}} + 2\rho^{\frac{p}{q}}\right] ds}\\
&\leq \frac{e^{(a_0 - a_{-})t} }{a_0}e^{\int_{\tau}^{t} b_0  \|J\|_{p} C_2(s) \left[ |\Omega|^{\frac{1}{q}} + 2\rho^{\frac{p}{q}}\right] ds}\| S -S_{0}\|_{p}.
\end{split}
\]
Thus, the result is immediate.
\eproof

\subsection{Upper semicontinuity of the pullback attractors}\label{upper-sc}

In this subsection  $\{\mathcal{A}_S(t);t\in\mathbb{R}\}$ denotes the pullback attractor for the process $T_S(\cdot,\cdot)$ in $X_p$, for $1 \leq p < \infty $.

Using Theorem \ref{tend}, we prove that
the family of pullback attractors $\{\mathcal{A}_S(t); t\in\mathbb{R}\}_{S\in \Sigma}$
is upper-semicontinuous at  $S_0 \in \Sigma$, i.e, we show that
\[
\lim_{t \to\infty} \operatorname{dist}_{H}(\mathcal{A}_S(t),
\mathcal{A}_{S_0}(t))=0,
\]
where $\operatorname{dist}_{H}(\cdot,\cdot)$ denotes the Hausdorff semi-distance.

\begin{theorem}\label{TheorUpper}
Under same hypotheses of Theorem  \ref{tend} the family of pullback attractors $\{\mathcal{A}_S(t); t\in\mathbb{R}\}_{S\in \Sigma}$ is upper-semicontinuous at $S_0 \in \Sigma$.
\end{theorem}
\proof
Note that, from Theorem \ref{Theor1}, it follows that
$$
\cup_{S\in \Sigma}\mathcal{A}_{S}(t) \subset B(0,R),
$$
where $R=R(t)=\dfrac{1}{a_{-}-k_1 b_0} [b_0 k_2(t) |\Omega|^{\frac{1}{p}} +p\| S\|_{p}]$.
Fixe $\varepsilon >0$ and $t\in \mathbb{R}$. Thus choose $\tau \in \mathbb{R}$, $\tau \leq t$, such that

$$
\operatorname{dist}_{H}(T_{S_0}(t,\tau)B(0,R),\mathcal{A}_{S_0}(t))<\dfrac{\varepsilon}{2}.
$$
Now, by Theorem \ref{tend}, it follows that there exists $\delta > 0$ such that, for $\|S-S_{0}\|_{p} < \delta$, we have

\[
\sup_{a_S\in \mathcal{A}_S(\tau)} \operatorname{dist}(T_S(t,\tau)a_S,T_{S_0}(t,\tau)a_S)<\dfrac{\varepsilon}{2}.
\]

Then, for $\|S-S_0\|_{p} < \delta$, using the invariance of the pullback attractors,  we obtain
\[
\begin{split}
&\operatorname{dist}_{H}(\mathcal{A}_S(t),\mathcal{A}_{S_0}(t))\\
&\leq \operatorname{dist}_{H}(T_{S}(t,\tau)\mathcal{A}_S(\tau),T_{S_0}(t,\tau)\mathcal{A}_S(\tau))+\operatorname{dist}_{H}(T_{S_0}(t,\tau)\mathcal{A}_S(\tau),T_{S_0}(t,\tau)\mathcal{A}_{S_0}(\tau))\\
&=\sup_{a_S\in \mathcal{A}_S(\tau)} \operatorname{dist}_{H}(T_S(t,\tau)a_S,T_{S_0}(t,\tau)a_S)+\operatorname{dist}_{H}(T_{S_0}(t,\tau)\mathcal{A}_S(\tau),\mathcal{A}_{S_{0}}(t))\\
&<\frac{\varepsilon}{2}+ \frac{\varepsilon}{2}=\varepsilon.
\end{split}
\] \eproof

\section{Discussions and Biological Interpretation} \label{conclude}


As we saw in the introduction, equation (\ref{equation}) generalizes the model studied by Amari (1977), 
which is already well known in the literature, because we consider that the rate in the intensity of neuronal potential is explicitly time dependent, 
while in Amari (1977), 
this rate was considered constant. We expect to have a more realistic model when compared to what happens in the brain, since this behavior is due to variations
of polarity inside the membrane, which is not necessarily constant. 
Furthermore, in Proposition \ref{globalexist} and Corollary \ref{Smoothness}, we are not considering that the synaptic connectivity function $J(x,y)$ is smooth, as occurs for example in Amari (1977, 1989), 
Bezerra at al. (2017) 
and Da Silva and Pereira (2015). 
For these results, we assume  $J\in L^1(\mathbb{R}^{N})$, leaving the model closer to real situation of mild autism, where occurs simple breaks in the synaptic connections. Thus, we hope that the results on global existence and smoothness of solutions, given in Proposition \ref{globalexist} and Corollary \ref{Smoothness} contribute to future research.

In Theorem 5.1 we show that the neuronal activity depends continuously on the sum of the external stimulus involved in the process. This reinforces the importance of appropriate continuous stimulation for a good neural activity, especially in individuals suffering from neurological disorder, as occurs in cases of cerebral paralysis and in some cases of autism.

Finally, we expect that the mathematical results presented in Theorem \ref{Theor1} and Theorem \ref{TheorUpper} contribute to other mathematical properties associated with the dynamics of this model and that other biological conclusions will be possible.

%
%
%



\end{document}